\documentclass[reqno,11pt]{amsart}

\newtheorem{theorem}{Theorem}[section]
\newtheorem{lemma}[theorem]{Lemma}
\newtheorem{corollary}[theorem]{Corollary}

\theoremstyle{definition}

\newtheorem{assumption}[theorem]{Assumption}
\newtheorem{example}[theorem]{Example}

\theoremstyle{remark}
\newtheorem{remark}[theorem]{Remark}

\newcommand{\mysection}[1]{\section{#1}
\setcounter{equation}{0}}

\newcommand{\WO}{\overset{\scriptscriptstyle0}%
{W}\,\!}

\newcommand{\bR}{\mathbb R}

\newcommand\cL{\mathcal{L}}

\newcommand{\nlimsup}{\operatornamewithlimits{\overline{lim}}}
\newcommand{\nliminf}{\operatornamewithlimits{\underline{lim}}}

\begin{document}
\title[Energy equality]
{On an energy equality in the theory of evolution
equations}
 
\author{N.V. Krylov}
\address{127 Vincent Hall, University of Minnesota, Minneapolis,
 MN, 55455}
 
\email{krylov@math.umn.edu}
 \keywords{Energy equality, spaces in duality}
 
\subjclass{35K90, 35R99}

\begin{abstract}
We give an elementary proof of
a basic ``energy'' equality in the theory
of evolution equations in the setting, which 
usually starts with introducing
 $V\subset H=H^{*}\subset V^{*}$.
 
\end{abstract}

\maketitle

\mysection{Main result, discussion, and   examples}

Let $V$ be a Banach space
over the field of complex numbers
 and let $H$ be a  Hilbert space 
over the field of complex numbers with 
 norms $\|\cdot\|_{V}$, $\|\cdot\|_{H}$ and scalar product
$(\cdot,\cdot)_{H}$, respectively. 
\begin{assumption}
                                   \label{assumption 11.15.1}
We have $V\subset H$ and for $v\in V$, $\|v\|_{V}\geq\|v\|_{H}$.
\end{assumption}

Let $V^{*}$ be the dual space to $V$
with pairing between them denoted by $\langle
v^{*},v\rangle$, $v\in V$, $v^{*}\in V^{*}$. Fix $T\in(0,\infty)$
and $p\in[1,\infty)$.
 
\begin{assumption}
                                   \label{assumption 11.15.2}
For $t\in[0,T]$,
 we are given a $V $-valued strongly (Lebesgue) measurable
function $v_{t}$ and a $V^{*}$-valued  
function $v^{*}_{t}$ such that 

(i) $\langle v^{*}_{t},v\rangle$
is a measurable function of $t\in[0,T]$
for any $v\in V$,

(ii)
there is a measurable function $f_{t}$ such that
$\|v^{*}_{t}\|_{V^{*}}\leq f_{t}$ for $t\in[0,T]$, 

(iii) there is a constant $N_{0}\in[1,\infty)$
such that 
\begin{equation}
                                                   \label{8.9.5}
 \int_{0}^{T}\big(\|v_{t}\|^{p}_{V}+f_{t}^{p'}  
\big)
\,dt\leq N_{0},  \quad p'=p/(p-1),\quad\text{if}\quad p>1,
\end{equation}
$$
 \int_{0}^{T} \|v_{r}\| _{V}\,dr+
f_{t}
\leq N_{0}\quad \text{(a.e.)}\,\,t\in[0,T]\quad\text{if}\quad p=1.
$$
\end{assumption}

\begin{assumption}
                                   \label{assumption 11.15.3}
 For any $v \in V $, for almost all $s,t\in[0,T]$, we have
\begin{equation}
                                          \label{10.26.1}
(v_{t},v )_{H}-( v_{s},v)_{H}
=\int_{s}^{t} \langle v_{r}^{*},v \rangle
\,dr.
\end{equation}

\end{assumption}
Our main goal is to present an elementary proof
of the following result.

\begin{theorem}
                                  \label{theorem 10.27.1}
Under the above assumptions there exists
an $H$-valued $H$-strongly continuous function $u_{t}$,
$t\in[0,T]$, such that $u_{t}=v_{t}$ for almost all
$t\in[0,T]$, and for all $t\in[0,T]$ and $v\in V$
\begin{equation}
                                          \label{10.26.2}
(u_{t},v )_{H}
=( u_{0},v)_{H}+\int_{0}^{t} \langle v_{s}^{*},v \rangle
\,ds.
\end{equation}
Furthermore, for all $t\in[0,T]$
\begin{equation}
                                          \label{10.27.1}
\|u_{t}\|^{2}_{H}
=\|u_{0}\|^{2}_{H}+2\int_{0}^{t}\Re \langle v_{s}^{*},u_{s}\rangle
\,ds.
\end{equation}

\end{theorem}

\begin{corollary}
                                       \label{corollary 11.6.1}
Under the above assumption,
suppose that $V$ is dense in $H$
(in the metric of $H$) and we are given
 $u\in H$ such that for any $v \in V $
$$
(v_{t},v )_{H}
=( u ,v)_{H}+\int_{0}^{t} \langle v_{r}^{*},v \rangle
\,dr
$$
for almost all $t\in[0,T]$. Then, of course,
\eqref{10.26.1} holds for almost all
$s,t\in[0,T]$ as long as $v\in V$,
there exists $u_{t}$ with the properties described 
in the theorem,
 and since $(u_{0},v)_{H}=(u,v)_{H}$ for any $v\in V$
and $V$ is dense in $H$,  we have $u_{0}=u$.
\end{corollary}

\begin{remark}
                                        \label{remark 11.6.1}
Obvious versions of Theorem \ref{theorem 10.27.1}
and Corollary \ref{corollary 11.6.1}
are true when $V$ and $H$ are spaces
 over the field of real numbers.
This will be easily seen from the proof of 
Theorem \ref{theorem 10.27.1}.
\end{remark}

The results like Theorem  \ref{theorem 10.27.1}
and Corollary \ref{corollary 11.6.1} are widely used
in the theory of nonlinear (quasi-linear)
 parabolic partial differential
equations satisfying a monotonicity condition.
These results serve as an intermediate tool
to achieve the goal of proving existence an uniqueness
of solutions. Probably because of that the author could
not find Theorem \ref{theorem 10.27.1}
and Corollary \ref{corollary 11.6.1} in the above
generality in the literature.

It seems the results like Theorem  \ref{theorem 10.27.1}
and Corollary \ref{corollary 11.6.1} first 
appeared in J.-L. Lions \cite{Li58} when $V$ was a Hilbert space.
However, the energy equality \eqref{10.27.1}
is not singled out because there were just no   point
in proving it. A different proof of the existence of 
$H$-continuous modifications using
mollification in $t$ is given in \cite{Li63}.
There the energy equality \eqref{10.27.1}
 is used as an obvious fact.

Then in Remark 2.1.2 on  page 156 of J.-L. Lions \cite{Li69} we read:
Let $V$ be a reflexive Banach space embedded into a Hilbert
space $H$, $V\subset H$, with continuous embedding,
and let $V$ be dense in $H$; after identifying $H$ with its dual
and denoting $V^{*}$ the dual of $V$, on can identify
$H$ with a subspace of $V^{*}$, so that $V\subset H\subset V^{*}$.
Then, if we are given a function $u\in L_{p}((0,T),V)$
such that $du/dt\in L_{p'}((0,T),V^{*})$, the function $u$
is (after being modified on a set of measure zero)
a strongly continuous $H$-valued function on $[0,T]$
and the mapping $u\to u(0)$ is onto $H$.

No further comments, references, or proofs are given,
and in the author's opinion, for an inexperienced reader even
to understand what exactly is claimed is a challenge.
Is it assumed that $p>2$ as in the title of the subsection 2.1.1?
Of course, the above statement in \cite{Li69}
is just a side remark and cannot be regarded
as a flaw of this  book remarkable in all respects.

Similar situation we have with Remark 2.7.9 on page 236 
of J.-L. Lions \cite{Li69}, where $V=\WO^{1}_{p}(\Omega)
\not\subset H=L_{2}(\Omega)$ 
(if $p\leq 2d/(d+2)$, $d=\text{dim}(\Omega)$),
$u\in L_{p}((0,T),V)$ and $du/dt\in
L_{p'}((0,T),V^{*})$. Still the claim is that
the function $u$
is (after being modified on a set of measure zero)
a strongly continuous $H$-valued function on $[0,T]$.
Again no comments or proofs are given and how the trivial
example of $u_{t}\equiv v\in V\setminus H$ is excluded
is not clear to an inexperienced reader.
We show a way to treat such situations in an almost
trivial generalization of Corollary \ref{corollary 11.6.1}
in Section \ref{section 11.19.2}.

The only place where the author could find the proof of
the results like Theorem  \ref{theorem 10.27.1}
and Corollary \ref{corollary 11.6.1} when $V$
is {\em not\/} a Hilbert space is the article by
F.E. Browder \cite{Br64}, where
$V$ is a closed subset of $W^{m}_{p}$, $p\in(1,\infty)$,
 hence a reflexive
space. He considers the closure of the operator
$d/dt$ as an operator from $X=L_{p}(\bR,V)$
to $X^{*}$. However, his setting is quite different from the one
based on $V\subset H=H^{*}\subset V^{*}$.

We decided to only concentrate on
  Theorem  \ref{theorem 10.27.1}
and Corollary \ref{corollary 11.6.1} and prove them
in the most general setting requiring
only basic knowledge of functional analysis
to understand the statements and their 
rather elementary proofs. We do not assume that $V$
is reflexive, or dense in $H$, or separable.
Our main tools
are some standard arguments from
measure theory, Doob's remarkable theorem
on approximation of Lebesgue integrals by
Riemannian sums, and the following two formulas
$$
(u_{2},v_{2})_{H}-(u_{1},v_{1})_{H}
=(u_{2}-u_{1},v_{1})_{H}+\overline{(v_{2}-v_{1},u_{2})_{H}}
$$
\begin{equation}
                                        \label{11.10.1}
\|u_{2}\|^{2}_{H}-\|u_{1}\|^{2}_{H}
=2\Re (u_{2}-u_{1},u_{1})_{H}+\|u_{2}-u_{1}\|^{2}_{H}
\geq 2\Re (u_{2}-u_{1},u_{1})_{H}.
\end{equation}

The proofs given here are adaptations of
the proof of It\^o's formula for Banach space-valued
stochastic processes from \cite{KR}, where  
$p\in(1,\infty)$. Unlike
   \cite{Li63} and \cite{Br64},
in the stochastic setting 
approximating $v_{t}$ with  functions
smoother
in time
does not lead to anything reasonable and one
could only use the tools mentioned above.
It is worth noting that a different proof of It\^o's formula
for Banach space-valued
stochastic processes is given earlier in
\cite{Pa75}. This proof is based on  
Remark 2.1.2 on  page 156 of J.-L. Lions \cite{Li69}.

We first give two examples of applications 
of Theorem \ref{theorem 10.27.1} in two border
situations, then prove this theorem
in Section \ref{section 11.19.1}. As pointed out above, in Section
\ref{section 11.19.2} we present a generalization of the main
theorem in the case that $V\not\in H$.

 Here is an example of application of
Theorem \ref{theorem 10.27.1} in which
$v^{*}_{t}$ is not strongly measurable
as a $V^{*}$-valued function. 

 \begin{example}
                                       \label{example 11.2.1}
For $q\in[1,\infty)$ let $\ell_{q}$ be the Banach space of
sequences $u=(u_{0},u_{1},....)$ of complex numbers
such that 
$$
\|u\|_{\ell_{q}}=\Big(\sum_{n=0}^{\infty}
|u_{n}|^{q}\Big)^{1/q}<\infty,
$$
and let $\ell_{\infty}$ be the Banach space of sequences
such that
$$
\|u\|_{\ell_{\infty}}=\sup_{n}|u_{n}|<\infty.
$$
Introduce $V=\ell_{1}$, then $V^{*}=\ell_{\infty}$
with pairing
$$
\langle v^{*},v\rangle=\sum_{n=0}^{\infty}
 v^{*}_{n}\bar v_{n}.
$$
Also introduce $H=\ell_{2}$ with the scalar product
$$
(u,v)=\sum_{n=0}^{\infty}u_{n}\bar v _{n}.
$$
Observe that $V\subset H$, $\|v\|_{V}\geq\|v\|_{H}$,
 and consider the function
$$
v_{t}=(v_{0}(t),v_{1}(t),...),\quad
v_{n}(t)= 2^{-n}\exp ( i2^{ n}\pi t), n=0,1,2,....
$$
Obviously, $v_{t}$ is a $V$-valued strongly continuous,
and hence measurable, function defined  for $t\geq0$
and 
$$
\int_{0}^{T}\|v_{t}\|_{V}^{2}\,dt<\infty.
$$
 Furthermore, 
for any $v\in V$
$$
(v_{t},v)_{H}=(v_{0},v)_{H}+\sum_{n=0}^{\infty} 
2^{-n}\big[\exp ( i2^{ n}\pi t)-1\big]
\bar v_{n}
$$
$$
=(v_{0},v)_{H}+i\pi \sum_{n=0}^{\infty} \int_{0}^{t}
\exp ( i2^{ n}\pi s)\bar v_{n}\,ds 
$$
$$
=(v_{0},v)_{H}+i\pi \int_{0}^{t}\sum_{n=0}^{\infty}
\exp ( i2^{ n}\pi s)\bar v_{n}\,ds=
(v_{0},v)_{H}+\int_{0}^{t}\langle v^{*}_{s},v \rangle\,ds,
$$
where
$$
v^{*}_{t}=(v^{*}_{0}(t),v^{*}_{1}(t),...),\quad
v^{*}_{n}(t)=i\pi\exp ( i2^{ n}\pi t), n=0,1,2,....
$$
Obviously,   $v^{*}_{t}$
is a  $V^{*}$-valued function, $\|v^{*}\|_{V^{*}}\leq
\pi=:f_{t}$,
$$
\int_{0}^{T}f^{2}_{t}\,dt<\infty,
$$
 and for any $v\in V$,
$\langle v^{*}_{t},v \rangle$ is  measurable. By Theorem
\ref{theorem 10.27.1} with $p=2$
$$
\|v_{t}\|_{H}^{2}=\sum_{n=0}^{\infty} 2^{-2n} =\|v_{0}\|^{2}
+2\Re\int_{0}^{t}\sum_{n=0}^{\infty}  2^{-n}i\pi\,ds=
\|v_{0}\|^{2}.
$$

This result is, of course, trivial in itself. What is remarkable
though is that $v^{*}_{t}$ is not a strongly measurable
$V^{*}$-valued function, and condition (i)
in Assumption \ref{assumption 11.15.2} is essential. 

To prove this fact, it suffices to show that
the set of values of  $v^{*}_{t}$ in $V^{*}$
as $t$ belongs to any set of full measure in $(0,T)$
is nonseparable.  We claim that a stronger assertion holds:
for any $s,t\in[0,1]$, such that $s\ne t$, we have
\begin{equation}
                                            \label{11.4.1}
\sup_{n=0,1,2,...}|\exp(i2^{n}\pi t)-\exp(i2^{n}\pi s)|
\geq  |1-\exp(i2\pi/3)| .
\end{equation}

 Assume that \eqref{11.4.1}
is false, set $r=t-s$,
and suppose without loss of generality
that $r>0$. Then, for any $n=0,1,2,...$,
there exists an integer $N_{n}$ and a number $\phi_{n}$
such that
$$
2^{n}\pi r= 2N_{n}\pi+\pi\phi_{n},\quad|\pi\phi_{n}|<
2\pi/3 .
$$
In short
$$
2^{n} r= 2N_{n}+ \phi_{n},\quad|\phi_{n}|< 2/3.
$$

Since $r>0$, $N_{0}=0$ and $\phi_{0}=r>0$.
Then for 
$$
n_{0}=\Big[\log_{2}\frac{2}{3\phi_{0} }\Big]+1
$$
where $[a]$ is the integer part of $a$, we have
 $
2^{n_{0}}r=2^{n_{0}}\phi_{0}$ and 
$$
2/3< 
2^{n_{0}}\phi_{0}\leq 4/3
$$
since $a-1<[a]\leq a$. It follows that 
in the representation $2^{n_{0}}r=2N_{n_{0}}+\phi_{n_{0}}$
($=2^{n_{0}}\phi_{0}$) with smallest possible $|\phi_{n_{0}}|$,
either $N_{n_{0}}=0$ and
$\phi_{n_{0}}=2^{n_{0}}\phi_{0}>2/3$ or $N_{0}=1$
and $\phi_{n_{0}}=2^{n_{0}}\phi_{0}-2\leq-2/3$.
This is the desired contradiction proving our claim.

\end{example}

Here is an example where again $V$ in non reflexive and $p=1$.

\begin{example}
 For $n=1,2,...$ set $c_{n}=2^{-n}$ and for
$x\in I:=(-\pi,\pi)$
 and $t\in[c_{n+1},c_{n})$ define
$$
k_{n}=[\ln_{2}n],\quad a_{n}=2^{n-k_{n}}, 
$$
$$
v_{t}(x)=(2^{n+1}t-1) c_{n}\sin (a_{n} x)
+(2 -t 2^{n+1}) c_{n+1}\sin (a_{n+1} x).
$$
Observe that $u_{t}$ is obtained by the linear interpolation
on $t\in[c_{n+1},c_{n})$ between the values at the end points,
which are $c_{n+1}\sin (a_{n+1} x)$ and $c_{n}\sin (a_{n} x)$.
The function $v_{t}(x)$ is defined in $  I  \times
(0,T)$, where $T=1/2$.

Note that
$$
\int_{-\pi}^{\pi}|D^{2}v_{c_{n}}|\,dx=c_{n}a_{n}^{2}
\int_{-\pi}^{\pi}|\sin (a_{n} x)|\,dx,
$$
which is of order $2^{n-2k_{n}}$ as $n\to\infty$. The series
$\sum 2^{-2k_{n}}$ converges and this implies that
$v_{\cdot}\in L_{1}((0,T),V)$, where $V=L_{2}( I )
\cap \WO^{2}_{1}( I )$ is a Banach space
with norm defined as the sum of norms in
$L_{2}(I)$ and $\WO^{2}_{1}(I)$. Obviously,
$V\subset H:=\WO^{1}_{2}(I)$.

Finally, the first derivative of $v_{t}(x)$ with respect to 
$t$, say $v^{*}_{t}(x)$, is a bounded function,
so that $v^{*}_{t}\in V^{*}$, and for $v\in C^{2}( I )$
vanishing at $\pm\pi$, by Fubini's theorem we have
$$
(v_{t},v)_{H}=\int_{-\pi}^{\pi}v_{t} (1-D^{2}v)\,dx
=\int_{0}^{t}\langle v^{*}_{s},v\rangle\,ds,
$$
where
$$
\langle v^{*}_{s},v\rangle=\int_{-\pi}^{\pi}
v^{*}_{s}(1-D^{2}v)\,dx
$$
extends to a bounded linear functional on $V$ with its norm
bounded on $(0,T)$. By Corollary \ref{corollary 11.6.1}
the function $v_{t}$ extends to $[0,T]$ as a strongly
continuous $H$-valued function equal to zero at $t=0$.

This is not a surprising result of course.
However, observe that it follows from 
\begin{equation}
                                                \label{11.19.1}
\int_{-\pi}^{\pi}|D^{2}v_{t}|\,dx
\geq (2^{n+1}t-1) \int_{-\pi}^{\pi}|D^{2}v_{c_{n}}|\,dx
-(2 -t 2^{n+1})\int_{-\pi}^{\pi}|D^{2}v_{c_{n+1}}|\,dx,
\end{equation}
valid for $t\in[c_{n+1},c_{n})$ , that there exists an 
absolute constant $\varepsilon\in(0,1/2)$ such that
if $t\in[c_{n}(1-\varepsilon), c_{n})$, then the
 right-hand side of \eqref{11.19.1} is bigger than
$\varepsilon 2^{n-2k_{n}}$. An easy consequence of this
is that $v_{\cdot}\not\in L_{p}((0,T),\WO^{2}_{q}( I ))$
for $p,q\geq1$ unless $p=q=1$.
\end{example}

\mysection{Proof of Theorem \protect\ref{theorem 10.27.1}}
                                          \label{section 11.19.1}

The proof of Theorem \ref{theorem 10.27.1}
is based on three lemmas.
The first one is just a reminder of part of basics
of integration theory.

\begin{lemma}
                                   \label{lemma 10.28.1}
Let $F$ be a Banach space and let $f_{t}$, $t\in[0,T]$,
be a strongly measurable $F$-valued function such that
for a $q\in[1,\infty)$ we have
$$
\int_{0}^{T}\|f_{t}\|^{q}_{F}\,dt<\infty.
$$
Set $f_{t}=0$ for $t<0$ and $t>T$. Then
$$
I:=
\lim_{a\to 0}\int_{0}^{T}\|f_{t}-f_{t+a}\|^{q}_{F}\,dt=0.
$$

\end{lemma}

Proof. By definition, for any $\varepsilon>0$
there exists an integer $n$, measurable sets $A_{k}\subset
[0,T]$, and $g_{k}\in F$, $k=1,...,n$, such that, for
$$
\phi_{t}:=\sum_{k=1}^{n}g_{k}I_{A_{k}}(t),
$$
we have
$$
\int_{0}^{T}\|f_{t}-\phi_{t}\|^{q}_{F}\,dt\leq\varepsilon.
$$
Also, it is known that for any measurable sets $A \subset
[0,1]$
$$
\lim_{a\to 0}\int_{0}^{T}\|I_{A }(t)-
I_{A }(t+a)\|^{q} \,dt=0.
$$
By also taking into account that 
$$
\int_{0}^{T}\|f_{t+a}-g_{t+a}\|^{q}_{F}\,dt
\leq\int_{0}^{T}\|f_{t }-g_{t }\|^{q}_{F}\,dt,
$$
we find
$$
\nlimsup_{a\to 0}\int_{0}^{T}\|f_{t}-f_{t+a}\|^{q}_{F}\,dt
\leq2\int_{0}^{T}\|f_{t}-\phi_{t}\|^{q}_{F}\,dt\leq 2\varepsilon,
$$
and the lemma is proved.

Denote $\kappa_{(-)} (n,x)=2^{-n}[2^{n}x]$,
$\kappa_{(+)} (n,x)=2^{-n}[2^{n}x]+2^{-n}$, where $[x]$
is the integer part of $x\in\bR$, $n=1,2,...$.

The following result, basically, belongs to
 J.~Doob.

\begin{lemma}
                                   \label{lemma 10.27.1}
Let the assumption of Lemma \ref{lemma 10.28.1} be satisfied.

Set $f_{t}=0$ for $t<0$ and $t>T$.
Then there exists
 a sequence of integers $n_{k}$, $k=1,2,...$,
 such that $n_{k}\to\infty$ and for almost any
$c\in(0,1)$ 
\begin{equation}
                                          \label{10.27.2}
\int_{0}^{T}\|f_{t}-f_{\kappa_{(\pm)}(n_{k},t+c)-c}\|^{q}_{F}\,dt
\to0
\end{equation}
as $k\to\infty$.
\end{lemma}

Proof. The function $\|f_{t}-f_{s}\|_{F}$ is a measurable
function of $(t,s)$. Hence, 
$\|f_{t}-f_{\kappa_{(\pm)}(n ,t+c)-c}\|_{F}$ are measurable as
 well for either sign + or -.
By Fubini's theorem
$$
\int_{0}^{1}\int_{0}^{T}\|f_{t}-f_{\kappa_{(\pm)}(n ,t+c)-c}\|_{F}^{q}
\,dtdc=
\int_{0}^{T} \Phi_{\pm n}(t )
\, dt,
$$
where
$$
\Phi_{\pm n}(t ):=\int_{0}^{1}\|f_{t}-
f_{t+\kappa_{(\pm)}(n ,t+c)-(t+c)}\|_{F}^{q}
\,dc=
\int_{t}^{t+1}\|f_{t}-f_{t+\kappa_{(\pm)}(n , c)-c}\|_{F}^{q}
\,dc
$$
$$
=
\int_{0}^{1}\|f_{t}-f_{t+\kappa_{(\pm)}(n , c)-c}\|_{F}^{q}
\,dc,
$$
where the last equality follows from the fact that
$\|f_{t}-f_{t+\kappa_{(\pm)}(n , c)-c}\|_{F}$ are periodic
functions of $c$ and one of their periods equals one.

Hence,
\begin{equation}
                                          \label{10.30.1}
\int_{0}^{1}\int_{0}^{T}\|f_{t}-
f_{\kappa_{(\pm)}(n ,t+c)-c}\|_{F}^{q}
\,dtdc=\int_{0}^{1}\Psi_{\pm n}( c)\,dc,
\end{equation}
where
$$
\Psi_{ \pm n}(c)=\int_{0}^{T}\|f_{t}
-f_{t+\kappa_{(\pm)}(n,c)-c}\|^{q}
_{F}\,dt.
$$
Observe that $\Psi_{\pm n}(c)\to0$
as $n\to\infty$ for any $c\in[0,1]$ by Lemma
\ref{lemma 10.28.1}.  
 Furthermore, obviously
$$
|\Psi_{\pm n}(c)|\leq 2^{q}\int_{0}^{T}\|f_{t} \|^{q}
_{F}\,dt.
$$
Hence in light of \eqref{10.30.1} and the 
dominated convergence theorem
$$
a(n,c):=\int_{0}^{T}\Big(\|f_{t}
-f_{\kappa_{(-)}(n ,t+c)-c}\|_{F}^{q}
+\|f_{t}-f_{\kappa_{(+)}(n ,t+c)-c}\|_{F}^{q}\Big)
\,dt\to0
$$
in $L_{1}([0,1])$. Then there is a subsequence
$n_{k}\to\infty$ such that $a(n_{k},c)\to0$
for almost any $c\in(0,1)$. This is exactly
what is asserted and the lemma is proved.
\begin{remark}
                                      \label{remark 11.12.1}

Observe that
$$
x\geq \kappa_{(-)} (n,x)\geq x-2^{-n},\quad
x\leq \kappa_{(+)} (n,x)\leq x+2^{-n},
$$
the functions $\kappa_{(\pm)}(t)(n,x)$
are right-continuous piece-wise constant and have jumps
at points $k2^{-n}$, $k=0,\pm1,...$. Also
it is useful to note that,
for a fixed $c\in(0,1)$, the graph of 
 the function $y=\kappa_{(-)} (n,x+c)-c$, $x\in\bR$,
is obtained from the graph of $y=\kappa_{(-)}(n,x)-x$, $x\in\bR$,
by sliding the latter appropriately along the diagonal
$y=x$ in the direction of lesser values of the coordinates.

Next, if we have a set $C\subset(0,T)$ of full measure, define $D=C\cup(-\infty,0)\cup(T,\infty)$ and observe that
for any $k=0,1,2,...$ and $n=1,2,...$ the point $k2^{-n}-c$
belongs to $D$ for almost any $c\in(0,1)$. It follows that
there exists a set $C_{0}$ of full measure in $(0,1)$ such that,
for any $c\in C_{0}$, all points  $k2^{-n}-c$ are in $D$ for all
$k=0,1,2,...$ and $n=1,2,...$, that is for any $c\in C_{0}$
$$
(0,T)\cap\{k2^{-n}-c:k=0,1,2,...,n=1,2,...\}\subset C
$$
Of course, we can take and fix $c_{0}\in C_{0}$ 
 so that \eqref{10.27.2} holds as well.

\end{remark}

 \begin{lemma}
                                            \label{lemma 10.10.1}
Under the assumptions of Theorem \ref{theorem 10.27.1}
there exists a set $C\subset(0,T)$ of full measure
such that the set $v_{C}:=\{v_{t}:t\in C\}$
 is separable in the metric of $V$ and  equation \eqref{10.26.1} holds for any
$v\in v_{C}$ and $s,t\in C$. Furthermore,
if the sequence $n_{k}$ is taken from Lemma \ref{lemma 10.27.1}
with $f_{t}=v_{t}$, $F=V$, and $q=p$, then
there exists $c_{0}\in(0,1)$ such that \eqref{10.27.2} holds and
all values of the functions
$$
\chi_{(\pm)}(k,t):=
\kappa_{(\pm)}(n_{k},t+c_{0})-c_{0}.
$$ 
which lie in $(0,T)$, belong to $C$
for any $k$, and the same is true for all points of jumps
of these functions $($coinciding with the values of 
$\chi_{(-)}(k,t))$.  
\end{lemma}

Proof. Since $v_{t}$ is a strongly measurable $V$-valued function,
there exists a set $A\subset(0,T)$ of full measure
such that $v_{A}=\{v_{t}:t\in A\}$ is a separable set
in $V$. Let $\{v^{i},i=1,2,...\}$ be a countable
subset of $v_{A}$ dense in $v_{A}$ in the metric of $V$.

Then observe that for each $v=v^{i}$ equation \eqref{10.26.1}
holds for almost all $(s,t)\subset A\times A$. Since we have
only countably many $v^{i}$'s, there exists
a set $B\in A\times A$ of full measure
such that, for any $v=v^{i}$,
equation \eqref{10.26.1}
holds for all $(s,t)\in B$. 
By Fubini's theorem and in light of the
fact that $A$ has full measure, there exists $s_{0}\in A$
such that  set $C=\{t\in (0,T):(s_{0},t)\in B\} $,
has full measure in $(0,T)$. Since $C\subset A$, $v_{C}
\subset v_{A}$, and $v_{C}$ is a separable subset of $V$.

Furthermore, by construction,
 for any $v\in \{v^{i},i=1,2,...\}$ and $t\in C$
we have  
$$
(v_{t},v )_{H}-( v_{s_{0}},v)_{H}
=\int_{s_{0}}^{t} \langle v_{r}^{*},v \rangle
\,dr.
$$
By subtracting such equalities we see that
\eqref{10.26.1} holds for any
$v\in \{v^{i},i=1,2,...\}$ and $s,t\in C$.
Now since 
$\{v^{i},i=1,2,...\}$ is dense in $v_{A}$ and, say
$$
\Big|\int_{s}^{t}\langle v^{*}_{r},v\rangle
\,dr-\int_{s}^{t}\langle v^{*}_{r},v^{i}\rangle
\,dr\Big|\leq \|v-v^{i}\|_{V}\int_{0}^{T}f_{s}\,ds,
$$
it follows that equation \eqref{10.26.1}
holds for all $s,t\in C$ and $v\in v_{A}=\{v_{t}:t\in A\}$
and, hence, for all $v\in v_{C}$ as well.
This proves the first statement of the lemma.
The second one follows directly from
Remark \ref{remark 11.12.1}.
The lemma is proved.

{\bf Proof of Theorem \ref{theorem 10.27.1}}.
Take $\chi_{(\pm)}(k,t)$ from Lemma \ref{lemma 10.10.1}
and introduce
$$
I_{k}=(0,T)\cap\{i2^{-n_{k}}-c_{0}:i=1,2,...\},
\quad \rho=\bigcup_{k=1}^{\infty}I_{k},
$$
so that $I_{k}$ is the set of values of $\chi_{(-)}(k,t)$, $t>0$,
 which lie in $(0,T)$, $I_{k},\rho\subset C$, $\chi_{(-)}(k,t)=t$,
$\chi_{(+)}(k,t)=t+2^{-n_{k}}$ for $t\in I_{k}$.

Next,
take $s,t\in I_{k}$ such that $s<t$ and $s,t$ are neighbors
in $I_{k}$ ($t-s=2^{-n_{k}}$). Then from Lemma \ref{lemma 10.10.1},
\eqref{10.26.1}
and the fact that $v_{s},v_{t}\in v_{C}$, because $s,t\in C$,
we get that
$$
\|v_{t}\|_{H}^{2}-\|v_{s}\|_{H}^{2}=
(v_{t}-v_{s},v_{s})_{H}+
\overline{(v_{t}-v_{s},v_{t})_{H}}
=\int_{s}^{t}\langle v^{*}_{r},v_{s}\rangle\,dr
+\int_{s}^{t}\overline{\langle v^{*}_{r},v_{t}\rangle}\,dr
$$
\begin{equation}
                                               \label{11.1.5}
=\int_{s}^{t}\langle v^{*}_{r},v_{\chi_{(-)}(k,r)}\rangle\,dr
+\int_{s}^{t}
\overline{\langle v^{*}_{r},v_{\chi_{(+)}(k,r)}\rangle}\,dr.
\end{equation}
By having in mind summing up such equalities as telescopic sums
we obtain that, if $s,t\in \rho$  are 
such that $s<t$,
 then there exists an integer $k_{0}$ for which 
$s,t\in I_{k_{0}}$
and 
\begin{equation}
                                         \label{10.31.3}
\|v_{t}\|_{H}^{2}-\|v_{s}\|_{H}^{2}=
 \int_{s}^{t}\langle v^{*}_{r},v_{\chi_{(-)}(k,r)}\rangle\,dr
+\int_{s}^{t}
\overline{\langle v^{*}_{r},v_{\chi_{(+)}(k,r)}\rangle}\,dr
\end{equation}
with $k=k_{0}$. Since obviously, $I_{k}$'s are nested
\eqref{10.31.3} also holds for all $k\geq k_{0}$.

 We now fix  $s,t\in\rho$ such that $s<t$ and let 
 $k\to\infty$ in \eqref{10.31.3}.
Observe that 
$$
\Big|\langle v^{*}_{r},v_{\chi_{(\pm)}(k,r)}\rangle
-\langle v^{*}_{r},v_{r}\rangle\Big|
\leq\|v_{\chi_{(\pm)}(k,r)}-v_{r}\|_{V}f_{r}
$$
and the right-hand side goes to zero in $L_{1}((0,T))$
as $k\to\infty$
owing to Assumption \ref{assumption 11.15.2}
 (iii) and H\"older's inequality
in case $p>1$. It follows that $\langle v^{*}_{r},v_{r} \rangle$
is a measurable function. Also
$$
\Big|
\int_{s}^{t}\langle v^{*}_{r},v_{\chi_{(\pm)}(k,r)}\rangle\,dr
-\int_{s}^{t}\langle v^{*}_{r},v_{r}\rangle\,dr
\Big|\leq\int_{s}^{t}\|v_{\chi_{(\pm)}(k,r)}-v_{r}\|_{V}f_{r}
\,dr
$$
and the right-hand side goes to zero as $k\to\infty$
by the above. As a result we see that, if
$s,t\in\rho$  are such that $s<t$, then
\begin{equation}
                                         \label{11.1.1}
\|v_{t}\|_{H}^{2}-\|v_{s}\|_{H}^{2}=
 2\int_{s}^{t}\Re\langle v^{*}_{r},v_{r}\rangle\,dr.
\end{equation}
It follows, in particular, that
$\|v_{t}\|_{H}^{2}$ is a uniformly continuous bounded
function on $\rho$.  

Next, we need a separable  subspace of $H$.
By Lemma \ref{lemma 10.10.1}
  the set $v_{C} $
 is separable in the metric of $V$. Then  $v_{C}$ is also
separable in the metric of $H$.
Its closure in $H$ is a separable Hilbert space, say $\hat H$,
and $v_{C}$
 is everywhere dense
in $\hat H$ in its metric inherited from $H$.

Now for $t\in\rho$ we define $ u_{t}=v_{t}$,
and if $t\in[0,T]\setminus\rho$, we take any sequence
$t_{m}\to t$ of {\em elements of\/} $\rho$ such that
$v_{t(m)}$ converges weakly in $\hat H$ to an element of $\hat H$
and call it $ u_{t}$
(we use   that $\rho\subset C $ and $v_{C}
\subset \hat H$).
 Observe that, in light of Lemma \ref{lemma 10.10.1},
for any $v \in v_{C} $ and $t\in [0,T]$, we have
\begin{equation}
                                          \label{11.1.4}
( u_{t},v )_{H}
=( u_{0},v)_{H}+\int_{0}^{t} \langle v_{s}^{*},v \rangle
\,ds.
\end{equation}
 Since $v_{C}$ is dense in $\hat H$, this of course implies
uniqueness in the definition of $ u_{t}$
for $t\in[0,T]\setminus\rho$ and also implies its weak continuity
on $[0,T]$.  

Next, take $t\in(0,T]$, $s\in\rho$, $s<t$,
 sufficiently large $k$, so that $s\in I_{k}$ and
there are points in $I_{k}$ which are less than $t$,
call $t_{k}$ the closest  element of $I_{k}$
to $t$ from the left (it may happen that $t_{k}=t$) and write
that owing to \eqref{11.1.4}, for any $v\in v_{C}$,
$$
( u_{t}- u_{t_{k}},v)_{H}
=\int_{t_{k}}^{t}\langle v_{r}^{*},v\rangle
\,dr.
$$
In light of\eqref{11.10.1}  
$$
\| u_{t_{1}}\|_{H}^{2}-\| u_{t_{2}}\|_{H}^{2} 
\geq 2\Re( u_{t_{1}}- u_{t_{2}}, u_{t_{2}})_{H},
$$
implying that  (recall that $u_{s}=v_{s}$ for $s\in\rho$)

$$
\| u_{t }\|_{H}^{2}-\| u_{t_{k}}\|_{H}^{2}
\geq2\int_{t_{k}}^{t }\Re\langle v^{*}_{r},v_{\chi_{(-)}(k,r)}
\rangle\,dr.
$$
Generally, for two neighboring points $t_{1},t_{2}\in I_{k}$
such that $t_{1}<t_{2}$
$$
\| u_{t_{2} }\|_{H}^{2}-\| u_{t_{1}}\|_{H}^{2}
\geq2\int_{t_{1}}^{t_{2}}\Re\langle v^{*}_{r},v_{\chi_{(-)}(k,r)}
\rangle\,dr,
$$
which,
 similarly to the manipulations after \eqref{11.1.5},
leads to
$$
\| u_{t }\|_{H}^{2}-\| u_{s}\|_{H}^{2}
\geq 2\int_{s}^{t}\Re\langle v^{*}_{r},v_{\chi_{(-)}(k,r)}  
\rangle\,dr.
$$
By what is said before this yields as $k\to\infty$ that
$$
\| u_{t }\|_{H}^{2}-\|  u_{s}\|_{H}^{2}
\geq 2\int_{s}^{t}\Re\langle v^{*}_{r},v_{r}
\rangle\,dr=\lim_{q\to t,q\in\rho}\|v_{q }\|_{H}^{2}-
\| v_{s}\|_{H}^{2}.
$$
However, since $ u_{t}$ is the weak limit
in $H$ of a sequence $v_{t(m)}$ such that $t(m)\in\rho$
and $t(m)\to t$, we have
$$
\| u_{t }\|_{H}^{2}\leq\nliminf_{m\to\infty}
 \|v_{t(m) }\|_{H}^{2}.
$$
It follows that
\begin{equation}
                                              \label{11.2.1}
\| u_{t }\|_{H}^{2}-\|  u_{s}\|_{H}^{2}
= 2\int_{s}^{t}\Re\langle v^{*}_{r},v_{r}
\rangle\,dr 
\end{equation}
as long as $t\in(0,T]$, $s\in\rho$, and $s<t$.

 We moved from $t\in(0,T]$ down to $s\in\rho$. Next we 
move   up to $s$ and we only need to do that by taking $t=0$.
So, fix $s\in\rho$, say $s\in I_{k}$, and denote
by $s_{k}$ the smallest element of $I_{k}$. Observe that
by \eqref{11.10.1}
$$
\| u_{0}\|^{2}_{H}-\| u_{s_{k}}\|^{2}_{H}
\geq 2\Re( u_{0}- u_{s_{k}},u_{s_{k}})_{H}= -2\int_{0}^{s_{k}}
\Re\langle v^{*},v_{\chi_{(+)}(k,r)}\rangle\,dr.
$$
Generally, for two neighboring points $s_{1},s_{2}\in I_{k}$
such that $s_{1}<s_{2}$
$$
\| u_{s_{1} }\|_{H}^{2}-\| u_{s_{2}}\|_{H}^{2}
\geq-2\int_{s_{1}}^{s_{2}}\Re
\langle v^{*}_{r},v_{\chi_{(+)}(k,r)}  
\rangle\,dr,
$$
which again by using telescoping sums leads to
$$
\| u_{0}\|^{2}_{H}-\| u_{s}\|^{2}_{H}
\geq -2\int_{0}^{s}\Re\langle v^{*}_{r},v_{\chi_{(+)}(k,r)}  
\rangle\,dr.
$$
By letting $k\to\infty$ we conclude
$$
\| u_{0}\|^{2}_{H}-\| u_{s}\|^{2}_{H}
\geq -2\int_{0}^{s}\Re\langle v^{*}_{r},v_{r}  
\rangle\,dr =\lim_{s\downarrow0,s\in\rho}
\|v_{s}\|^{2}_{H}-\| u_{s}\|^{2}_{H}.
$$
This shows that
$$
\| u_{0}\|^{2}_{H}-\| u_{s}\|^{2}_{H}
= -2\int_{0}^{s}\Re\langle v^{*}_{r},v_{r}  
\rangle\,dr 
$$
and along with \eqref{11.2.1} yields that, for all $t\in[0,T]$,
\begin{equation}
                                              \label{11.2.2}
\| u_{t }\|_{H}^{2}=\|  u_{0}\|_{H}^{2}
+ 2\int_{0}^{t}\Re\langle v^{*}_{r},v_{r}  
\rangle\,dr .
\end{equation}
Hence, $\| u_{t }\|^{2}_{H}$ is continuous
on $[0,T]$, and, since $ u_{t}$ is weakly continuous,
it is  strongly
continuous on $[0,T]$ as an $\hat H$-valued function
and as an $H$-valued function as well.

Next, we claim  that 
$ u_{t}=v_{t}$ (a.e.) on $(0,T)$,
and, in particular, that 
one can replace
$v_{r}$ in \eqref{11.2.2} with $ u_{r}$.

By using Lemma \ref{lemma 10.10.1} and
 setting $s\downarrow0$
along $\rho$ in equation \eqref{10.26.1},
we get that for all $t\in C$
equation \eqref{11.1.4} holds if we replace
$u_{t}$ in the left-hand side with $v_{t}$.
Hence,
\begin{equation}
                                              \label{11.2.3}
( u_{t},v)_{H}=
( v_{t},v)_{H}
\end{equation}   
for all $t\in C$ and $v\in v_{C}$. 
Furthermore, since $v_{C}\in\hat H$, we have
$v_{t}\in\hat H$ if $t\in C$, and $u_{t}\in\hat H$ for all $t\in[0,T]$.
In addition $v_{C}$
is dense in  $\hat H$ in the metric of $\hat H$.
It follows that, for $t\in C$, \eqref{11.2.3}
holds for any $v\in \hat H$ meaning that $u_{t}=v_{t}$
for $t\in C$, that is almost everywhere  on $(0,T)$.
This proves our claim.

Finally, we deal with \eqref{10.26.2}. Since
$u_{t}=v_{t}$ almost everywhere, for any $v \in V $, we have
\begin{equation}
                                       \label{11.6.3}
(u_{t},v )_{H}-( u_{s},v)_{H}
=\int_{s}^{t} \langle v_{r}^{*},v  \rangle
\,dr
\end{equation}
for almost all $s,t\in[0,T]$. Here,
for any $v \in V $, the left-hand side
is a continuous functions of $t,s$ owing to the
strong continuity of $u_{t}$ in   $H$.
It follows that \eqref{11.6.3} holds true
for any $s,t\in[0,T]$ and $v \in V $, and
this brings the proof of the theorem to an end.

\mysection{A generalization}
                                            \label{section 11.19.2}

Here we keep Assumption \ref{assumption 11.15.2}
and replace Assumptions \ref{assumption 11.15.1}
and \ref{assumption 11.15.3} with the following ones.
If $F,G$ are Banach spaces, then by $\cL(F,G)$ we denote the space 
of bounded linear operators from $F$ to $G$.

\begin{assumption}
                                 \label{assumption 11.15.4}
We have that $\hat V:=V\cap H$ is dense in $H$ in the metric of $H$.
\end{assumption}
We introduce a norm on $\hat V$ by $\|\cdot\|_{\hat V}
=\|\cdot\|_{V }+\|\cdot\|_{  H}$ which makes
$V\cap H$ a Banach space.

\begin{assumption}
                                 \label{assumption 11.15.5}
There are   sequences $M_{k}$, $M'_{k}$, $k=1,2,...$, 
of linear operators  $M_{k}\in \cL(V,\hat V)$
(so to speak, $M_{k}$ are mollifying operators)
and $M'_{k}\in \cL( \hat V,V)$ (kind of adjoint to $M_{k}$)
such that
 
(i)  The norms of operators $M'_{k}M_{k}$ as elements
of $\cL(V,V)$ are bounded by a constant independent of $k$.
 
(ii) As $k\to\infty$, $M_{k}u\to u$ in $H$
if $u\in H$, and $M'_{k}M_{k}v\to v$ in $V$
if $v\in V$;

(iii) If $v\in V$, $u\in H$, and $M_{k}v\to u$ in $H$
as $k\to\infty$, then $v=u$.
\end{assumption}
\begin{assumption}
                                 \label{assumption 11.15.6}
There is a $w\in H$ such that,
for any $k =1,2,...$ and $v \in V\cap H$,
for almost all $ t\in[0,T]$, we have
\begin{equation}
                                          \label{11.15.1}
(M_{k}v_{t}, v )_{H}=(M_{k}w ,v)_{H}
=\int_{s}^{t} \langle v_{r}^{*},M'_{k} v \rangle
\,dr.
\end{equation}
 
\end{assumption}
 
\begin{theorem}
                                         \label{theorem 11.15.1}
Under the above assumptions
there exists
an $H$-valued $H$-strongly continuous function $u_{t}$,
$t\in[0,T]$, such that $u_{t}=v_{t}$ for almost all
$t\in[0,T]$, $u_{0}=w$, and for all $t\in[0,T]$ 
\begin{equation}
                                          \label{11.15.2}
\|u_{t}\|^{2}_{H}
=\|u_{0}\|^{2}_{H}+2\int_{0}^{t}\Re \langle v_{s}^{*},u_{s}\rangle
\,ds.
\end{equation}
\end{theorem}
 
Proof. Fix $k=1,2,..$, and on $v\in\hat V$ define
   $\langle\hat v^{*}_{t},v\rangle:=\langle v^{*}_{t},M'_{k}v
\rangle$. Observe that $\hat v^{*}_{t}\in \hat V^{*}$.
By applying Corollary \ref{corollary 11.6.1}
to $\hat V$, $M_{k}v_{t}$, and $\hat v^{*}_{t}$ in place of $V$,
$v_{t}$, and 
$v^{*}_{t}$, respectively we conclude that there exists
an $H$-valued $H$-strongly continuous function $u^{k}_{t}$,
$t\in[0,T]$, such that $u^{k}_{t}=M_{k}v_{t} $ for almost all
$t\in[0,T]$, $u^{k}_{0}=M_{k}w$, and for all $t\in[0,T]$ 
\begin{equation}
                                          \label{11.15.3}
\|u^{k}_{t}\|^{2}_{H}
=\|M_{k}w\|^{2}_{H}
+2\int_{0}^{t}\Re \langle v_{s}^{*},M'_{k}M _{k}v_{s}\rangle
\,ds.
\end{equation}

By applying the same argument to $M_{k}v_{t}
-M_{r}v_{t}$ we see that  for all $t\in[0,T]$ 
\begin{equation}
                                          \label{11.15.4}
\|u^{k}_{t}-u^{r}_{t}\|^{2}_{H}
=\|M_{k}w-M_{r}w\|^{2}_{H}
+2\int_{0}^{t}\Re \langle v_{s}^{*},
M'_{k}M _{k}v_{s}-M'_{r}M _{r}v_{s}\rangle
\,ds.
\end{equation}
It follows that uniformly on $[0,T]$ functions
$u^{k}_{t}$ converge in $H$ to an $H$-valued continuous
function $u_{t}$.

On the set of full measure in $(0,T)$ we have
$u^{k}_{t}=M_{k}v_{t}\to u_{t}$ in $H$. By assumption,
$v_{t}=u_{t}\in H$ on this set.
Passing to the limit in \eqref{11.15.3} on the basis
of the dominated convergence theorem presents
no difficulty and the theorem is proved.


\begin{thebibliography}{mm}

\bibitem{Br64} F.E. Browder, 
{\em Strongly non-linear parabolic boundary value
problems\/}, American J. of Math, Vol. 86 (1964),
No. 2, 339--357.
 
 \bibitem{KR} N.V. Krylov and B.L. Rozovsky,  {\em  
Stochastic evolution
equations},  ``Itogy nauki i
tekhniki'',  Vol. 14, VINITI, Moscow, 1979, 71-146 in Russian;
English translation in
 J. Soviet Math., Vol. 16 (1981), No. 4, 1233-1277. 

\bibitem{Li58} J.-L. Lions, {\em Espaces
interm\'ediaires entre espaces hilbertiens et applications\/}, 
Bull. Math. R.P.R. Bucarest, 2 (1958), 419--432.



\bibitem{Li63} J.-L. Lions, {\em\'Equations diff\'erentielles
op\'erationnelles
dans les
espaces de Hilbert\/}, 
 Equazioni differenziale astratte, Varenna, Italy, 1963,
Lectures given at Centro Internazionale 
Matematico Estivo [held in Varenna (Como),
Italy,
May 30-June 8 1963] , pp. 47- 122,
  Springer, 2011,


 \bibitem{Li69} J.-L. Lions,
``Quelques m\'ethodes de r\'esolution des probl\`emes aux limites
non lin\'eaires'', Dunod, Gauthier Villars, 1969.

\bibitem{Pa75} E.~Pardoux, {\em \'Equations aux d\'eriv\'ees
partielles  stochastiques non lin\'eaires monotones\/}, Ph.
D. Thesis, Universit\'e de Paris  Sud, Orsay, 1975,
http://www.cmi.univ-mrs.fr/~pardoux/Pardoux\_these.pdf

 

\end{thebibliography}
\end{document}